\documentclass[11pt, a4paper]{article}
\usepackage{amsmath}
\usepackage{fancyhdr}
\usepackage[misc]{ifsym}
\usepackage{color}
\usepackage{amsmath}
\usepackage{amsfonts}
\usepackage{dsfont}
\usepackage{mathrsfs}
\usepackage{bbm}
\usepackage{latexsym}
\usepackage{graphicx}
\usepackage{epstopdf}
\usepackage{amssymb}
\usepackage{indentfirst}
\usepackage{subfigure}
\usepackage{subcaption}
\usepackage{booktabs}
\usepackage{authblk}
\usepackage{cite}
\usepackage{float}
\usepackage[numbers,sort&compress]{natbib}
\usepackage{caption}
\usepackage{booktabs}
\usepackage{threeparttable}
\usepackage{multicol}
\usepackage{multirow}
\usepackage{geometry}
\newtheorem{remark}{Remark}
\newtheorem{theorem}{Theorem}[section]

\newtheorem{definition}[theorem]{Definition}

\numberwithin{equation}{section}

\def\v{\mathbf}

\begin{document}

\title{More than one Author with different Affiliations}
\author{Xuyang Na\\
\texttt{nxy.compmath@qq.com}
}

\title{An optimal two-side Robin-Robin domain decomposition method for H(div)-elliptic problem}\date{}
\maketitle

{\bf{Abstract}:} In this paper, we develop a new two-side Robin-Robin domain decomposition method for H(div)-elliptic problem. In \cite{zeng2015parallel}, Chen applied a Robin-Robin domain decomposition iterative method on H(div)-elliptic problem and proved the convergence rate is $1-O(h)$, where $h$ is the mesh size. We improve the algorithm by giving appropriate constraints on the interfaces of subdomains. Numerical results show that the convergence rate of the new algorithm only depends on $H/h$, where $H$ is diameter of subdomains. Besides, an algebraic system of Robin boundary conditions is derived from the iterative method. We solve it by MINRES and get asymptotically stable iteration numbers as well.

{\bf{Keywords}:} \hspace*{2pt}Robin-Robin domain decomposition method,\ \ H(div)-elliptic problem,\ \ Raviart-Thomas finite element.

\section{Introduction}

Robin-Robin nonoverlapping domain decomposition method, based on using Robin boundary conditions as transmission data on the interfaces of subdomains, has become one of the most important algorithms for solving linear systems derived from partial differential equations over the past few decades. The idea of employing Robin-Robin coupling conditions in domain decomposition was first proposed by P.L.Lions in \cite{lions1990schwarz}. The convergence without any rate is shown in \cite{lions1990schwarz} and \cite{quarteroni1999domain}. Later, Douglas, Huang\cite{douglas1997accelerated, douglas1998accelerated} and Guo, Hou\cite{guo2003generalizations} proved that the convergence rate could be $1-O(h)$, where $h$ is the mesh size. Gander\cite{gander2000optimized} pointed that the optimal choice of Robin parameter is $\gamma = O(h^{-1/2})$ by using Fourier analysis and the convergence rate $1-O(h^{1/2})$ could be achieved. Later, Xu and Qin\cite{xu2010spectral} gave a rigorous proof on this result and proved that the convergence rate is asymptotically optimal. At the same time, some researchers considered to extend the method to the case of many subdomains and solve various problems. Qin and Xu\cite{qin2006parallel} developed a Robin-Robin domain decomposition method with many subdomains of Poisson equation and proved that the convergence rate is $1-C^Nh^{1/2}H^{-1/2}$, where $H$ is the diameter of subdomains and $N$ is wind number. Bendali and Boubendir\cite{bendali2006non} established a unified framework for proving the convergence rate in the case of many subdomains. In subsequent work\cite{benhassine2011non}, Benhassine and Bendali extended the method in \cite{bendali2006non} to Stokes equations with continuous pressure. They got a convergence rate depend on the mesh size $h$. Later, Zeng, Chen and Li\cite{zeng2015parallel} solved H(div)-elliptic problem by this method without coarse problem and they proved the convergence is $1-O(h)$. Besides, as the format of Robin boundary condition, Robin-Robin domain decomposition method and its extensions are widely employed for solving Helmholtz equations\cite{boubendir2018non, liu2019optimized, claeys2022robust}. However, all the literature above could not provide an optimal convergence rate.\\
\indent
In this paper, we introduce a new Robin-Robin domain decomposition method with many subdomains for H(div)-elliptic problem. It is shown by numerical experiments that the convergence rate of our algorithm is optimal, that is, the convergence rate depends solely on the ratio between the diameter of subdomains $H$ and mesh size $h$. The outline is as follows: In section 2, we introduce the model problem. Section 3 describes how we get the new algorithm from present algorithms. Besides, we implement our new algorithm in two ways. In section 3, some numerical experiments are carried out to show the optimality of the new algorithm.

\section{Model problem}
We consider the boundary value problem for a vector field
\begin{equation}\label{sec1_hdiv_pro}
  \begin{cases}
  \begin{array}{rll}
    -\v{grad}(\text{div}\v{u})+\beta\v{u} & = \v{f}\quad &\text{in\ }\Omega, \\
    \v{u}\cdot\v{n} & = 0\quad &\text{on\ }\partial\Omega, 
  \end{array}
  \end{cases}
\end{equation}
where $\Omega$ is a bounded polygonal domain in $R^2$, $\v{n}$ is its outwards normal direction, $\beta$ is a positive number and $\v{f}\in (L^2(\Omega))^2$.\\
\indent
Let
\begin{equation*}
  H(\text{div};\Omega) = \{\v{u}\in (L^2(\Omega))^2:\text{div}\v{u}\in L^2(\Omega)\},
\end{equation*}
be equipped with the following inner product and the graph norm
\begin{gather*}
  (\v{u},\v{v})_{div} = \int_{\Omega}\text{div}\v{u}\text{div}\v{v}\text{d}x+\int_{\Omega}\v{u}\cdot\v{v}\text{d}x,\\
  \Vert \v{u}\Vert_{div}^2 = (\v{u},\v{u})_{div}.
\end{gather*}
Let $H_0(\text{div};\Omega)$ denote the subspace of $H(\text{div};\Omega)$ with vanishing normal components on the boundary $\partial\Omega$, i.e.
\begin{equation*}
  H_0(\text{div};\Omega) = \{\v{u}\in H(\text{div};\Omega):\v{u}\cdot\v{n} = 0\}.
\end{equation*}
The weak formulation of (\ref{sec1_hdiv_pro}) is to find $\v{u}\in H_0(\text{div};\Omega)$ such that
\begin{equation}\label{sec1_stokes_weak}
    a(\v{u},\v{v}) = (\v{f},\v{v})\quad \forall \v{v}\in H_0(\text{div};\Omega)
\end{equation}
where the bilinear forms are defined as follows:
\begin{gather*}
  a(\v{u},\v{v}) = \int_{\Omega}\text{div}\v{u}\text{div}\v{v}\text{d}x+\beta\int_{\Omega}\v{u}\cdot\v{v}\text{d}x,\quad\v{u},\v{v}\in H(\text{div};\Omega),\\
  (\v{f},\v{v}) = \int_{\Omega}\v{f}\cdot\v{v}\text{d}x,\quad\v{v}\in H(\text{div};\Omega).
\end{gather*}
Here, we notice that the energy norm $\Vert\cdot\Vert_a$ defined by $\Vert\cdot\Vert_a^2 = a(\cdot,\cdot)$ is equivalent to the graph norm as $\beta$ is a positive constant.

\section{Two-side Robin-Robin domain decomposition method for $H(\text{div})$ problem}

\subsection{Geometry settings and finite element spaces.}
Let $\mathcal{T}_h$ be a quasi-uniform and shape regular triangulation of the domain $\Omega$ with mesh size $h$ and $\mathcal{E}_h$ is the set of edges of triangular elements in $\mathcal{T}_h$. For every edge $e\in \mathcal{E}_h$, we fix a direction, given by a unit vector $\v{n}_e$, normal to $e$. The length of the edge $e$ is denoted by $\vert e\vert$.\\
\indent
Assume that $\Omega$ could be decomposed into some connected subdomains $\Omega_i(i = 1,\dots,N)$ satisfying
\begin{gather*}
  \Omega_i\cap\Omega_j = \phi,\quad i\neq j,\\
  \Omega = \bigcup_{i = 1}^{N}\overline{\Omega}_i,
\end{gather*} 
and the partition is compatible with the triangulation $\mathcal{T}_h$, that is, each element of $\mathcal{T}_h$ is contained in one of the subdomain $\Omega_i$. The diameter of $\Omega_i$ is $H_i$ and $H$ is the maximum of the diameters of the subdomains.\\
\indent
The edges of subdomains are denoted by
\begin{equation*}
  \Gamma_{ij} = \Gamma_{ji} = \partial\Omega_i\cap\partial\Omega_j,\quad i\neq j, \quad \vert \Gamma_{ij}\vert>0,
\end{equation*}
where $\vert \Gamma_{ij}\vert$ is the measure of $\Gamma_{ij}$. For every $\Gamma_{ij}$, we fix the direction of its normal vector and denote it by $\v{n}_{ij}$. Moreover, for convenience, we suppose that the fixed normal vector $\v{n}_e$ of each $e\in\mathcal{E}_h$ has the same direction with $\v{n}_{ij}$ if $e\subset\Gamma_{ij}$.
We denote $\mathcal{E}_H$ the set of edges of the subdomains and let $\Gamma$ be the union of the edges in $\mathcal{E}_H$, in other words, the subdomain boundaries do not belong to $\partial\Omega$:
\begin{equation*}
  \Gamma = \bigcup_{i = 1}^N\partial\Omega_i\backslash\partial\Omega.
\end{equation*}
For every subdomain $\Omega_i$, we denote $\mathcal{I}_i$ the set of indices $j$ such that $\Gamma_{ij}$ is an edge of $\Omega_i$:
\begin{equation*}
  \mathcal{I}_i = \{j| \Gamma_{ij}\subset\partial\Omega_i,\Gamma_{ij}\in\mathcal{E}_H\}.
\end{equation*}

In this article, we employ the lowest order Raviart-Thomas finite element space $\widehat{W}$ defined by
\begin{equation*}
  \widehat{W} = \{\v{u}\in H(\text{div};\Omega):\v{u}_{|K} \in \mathcal{RT}(K),K\in\mathcal{T}_h\}
\end{equation*}
where $\mathcal{RT}(K)$ is given by
\begin{equation*}
  \mathcal{RT}(K) := \v{a}+b\v{x},
\end{equation*}
and $K$ is an arbitrary element of $\mathcal{T}_h$. The degrees of freedom are defined by the average values of the normal components over the edges in $\mathcal{E}_h$, i.e.
\begin{equation*}
  \lambda_e(\v{u}) = \frac{1}{\vert e\vert}\int_e\v{u}\cdot\v{n}_e\text{d}s = \v{u}\cdot\v{n}_e,\quad e\in\mathcal{E}_h.
\end{equation*}
Then the finite element approximation of (\ref{sec1_hdiv_pro}) is to find $\v{u}\in \widehat{W}_0$ such that
\begin{equation}\label{hdiv_discrete}
  a(\v{u},\v{v}) = (\v{f},\v{v})\quad \forall\v{v}\in \widehat{W}_0
\end{equation}
where $\widehat{W}_0 := \widehat{W}\cap H_0(\text{div};\Omega)$ is the subspace of $\widehat{W}$ with vanishing normal component on the boundary. \\
\indent
Based on the subdomain partition, local function spaces could be defined as follows:
\begin{gather*}
  H_{\star}(\text{div};\Omega_i) = \{\v{u}_i\in H(\text{div};\Omega_i)| \v{u}_i\cdot\v{n} = 0 \text{ on } \partial\Omega_i\cap\partial\Omega\},\\
  W_i = {\widehat{W}_0}|_{\Omega_i} = \{ \v{u}_i\in H_{\star}(\text{div};\Omega_i)|{\v{u}_i}_{|_K}\in \mathcal{RT}(K),K\in\mathcal{T}_h,K\subset\Omega_i\}.
\end{gather*}
We next introduce some trace spaces $V_i$ consisting of normal components on the boundaries of the subdomains. A scalar function $u_i$, defined on $\partial\Omega_i\backslash\partial\Omega$, belongs to $V_i$ if and only if there exists $\v{u}_i\in W_i$ such that, for each edge $\Gamma_{ij}$ of $\Omega_i$,
\begin{equation*}
  {u_i}_{|_{\Gamma_{ij}}} = \v{u}_i\cdot\v{n}_{ij},\quad \Gamma_{ij}\in\mathcal{E}_H, j\in\mathcal{I}_i.
\end{equation*}
In other words, functions in $V_i$ are piecewise constant along the edges of $\Omega_i$.

\subsection{New Robin-Robin domain decomposition method}
Let $\v{u}\in H_0(\text{div},\Omega)$ be the solution of Equation (\ref{sec1_hdiv_pro}) and $\v{u}_i(i = 1,\dots,n)$ be the restriction of $\v{u}$ on subdomain $\Omega_i$. Then $\v{u}_i(i = 1,\dots,n)$ satisfies the following subproblem,
\begin{equation}\label{subpro}
  \begin{cases}
  \begin{array}{rll}
    -\v{grad}(\text{div}\v{u}_i)+\beta\v{u}_i & = \v{f}_i\quad &\text{in\ }\Omega_i, \\
    \v{u}_i\cdot\v{n} & = 0\quad &\text{on\ }\partial\Omega_i\cap\partial\Omega,\\
    \v{u}_i\cdot\v{n}_{ij} &= \v{u}_j\cdot\v{n}_{ij}\quad &\text{on }\Gamma_{ij},\ j\in\mathcal{I}_i,  \\
    \text{div}\v{u}_i &= \text{div}\v{u}_j\quad &\text{on }\Gamma_{ij},\ j\in\mathcal{I}_i.
  \end{array}
  \end{cases}
\end{equation}
where $\v{n}$ is the outer normal vector of domain $\Omega$.\\
\indent
Combine the third and forth equations of (\ref{subpro}), we may get two Robin-type conditions $g_{ij},g_{ji}\in V_{ij}$ as follows
\begin{align}
  \label{robinbd1} g_{ij}:= \gamma_{ij}\v{u}_i\cdot\v{n}_{ij}+\text{div}\v{u}_i & = \gamma_{ij}\v{u}_j\cdot\v{n}_{ij}+\text{div}\v{u}_j \\
  \label{robinbd2}g_{ji}:= \gamma_{ji}\v{u}_i\cdot\v{n}_{ij}-\text{div}\v{u}_i & = \gamma_{ji}\v{u}_j\cdot\v{n}_{ij}-\text{div}\v{u}_j
\end{align}
where $\gamma_{ij},\gamma_{ji}$ are arbitrary positive numbers.\\
\indent
Using the Robin-type conditions (\ref{robinbd1}) and (\ref{robinbd2}), we could get the following subproblem which is equivalent to (\ref{subpro}),
\begin{equation}\label{subrobin}
  \begin{cases}
  \begin{array}{rll}
    -\v{grad}(\text{div}\v{u}_i)+\beta\v{u}_i & = \v{f}_i\quad &\text{in\ }\Omega_i, \\
    \v{u}_i\cdot\v{n} & = 0\quad &\text{on\ }\partial\Omega_i\cap\partial\Omega,\\
    \gamma_{ij}\v{u}_i\cdot\v{n}_{ij}+\text{div}\v{u}_i & = g_{ij}\quad &\text{on }\Gamma_{ij},\ j\in\mathcal{I}_i,
  \end{array}
  \end{cases}
\end{equation}
together with the relation of Robin-type condition on the artificial interface $\Gamma_{ij}$
\begin{equation}\label{trans}
  g_{ji} = (\gamma_{ij}+\gamma_{ji})\v{u}_i\cdot\v{n}-g_{ij}.
\end{equation}
The weak form of problem (\ref{subrobin}) is to find $\v{u}_i\in W_i$ such that
\begin{equation}\label{subproWeak}
  a(\v{u}_i,\v{v}_i)+\sum_{j\in\mathcal{I}_i}\gamma_{ij}\langle \v{u}_i\cdot\v{n}_{ij},\v{v}_i\cdot\v{n}_{ij}\rangle = (\v{f},\v{v}_i)+\sum_{j\in\mathcal{I}_i}\langle g_{ij},\v{v}_i\cdot\v{n}_{ij}\rangle\quad \forall \v{v}_i\in W_i,
\end{equation}
where
\begin{equation*}
  \langle p,q\rangle = \int_{\Gamma}pq\text{d}s\quad p,q\in V = \Pi_{i = 1}^N V_i.
\end{equation*}

Now we set all the $\gamma_{ij}$ to the same value $\gamma$ and decompose $\v{u}_i$ into two parts, $\v{u}_{i,I}$ and $\v{u}_{i,\Delta}$. Then the problem (\ref{subproWeak}) could be written in matrix form as follows,
\begin{equation}\label{subproMat}
  \begin{pmatrix}
    A^i_{II} & A^i_{I\Delta} \\
    A^i_{\Delta I} & A^i_{\Delta\Delta}+\gamma M_i
  \end{pmatrix}
  \begin{pmatrix}
    u_{i,I} \\
    u_{i,\Delta} 
  \end{pmatrix} = 
  \begin{pmatrix}
    f_{i,I} \\
    f_{i,\Delta}+M_ig_i 
  \end{pmatrix}.
\end{equation}
Here the subscripts $I,\Delta$ denote the interior degrees of freedom and those on $\partial\Omega_i\backslash\partial\Omega$ respectively and the elements in matrix $A, M$ are obtained by computing the integrations $a(\cdot,\cdot),\langle\cdot,\cdot\rangle$ with corresponding basis functions in $W_i, V_i$. Besides, as the vector functions in $W_i$ are uniquely determined by the degrees of freedom defined on the edges in $\mathcal{E}_h$, we denote a vector function $\v{u}_i$ in $W_i$ by its corresponding column vector of degrees of freedom, named $u_i$, throughout the paper.\\
\indent
Gather all the subproblems on subdomains, we get
\begin{equation}\label{matAlgo1}
  \begin{pmatrix}
    A_{II} & A_{I\Delta} \\
    A_{\Delta I} & A_{\Delta\Delta}+\gamma M 
  \end{pmatrix}
  \begin{pmatrix}
    u_I \\
    u_{\Delta} 
  \end{pmatrix} = 
  \begin{pmatrix}
    f_I \\
    f_{\Delta}+Mg 
  \end{pmatrix}.
\end{equation}
In \cite{zeng2015parallel}, Chen introduce an iterative algorithm for H(div)-elliptic problem by solving (\ref{matAlgo1}) and updating the Robin boundary condition $g$ with the equality (\ref{trans}) in every iteration step. In their algorithm, the function space $W_0$ is decomposed into several local spaces $W_i$ and the solution is searched in the product space $W = \Pi_{i = 1}^N W_i$. However, the convergence rate of their algorithm deteriorates as the mesh size $h$ decreases. \\
\indent
To improve the performance, we develop a new iterative algorithm by searching the solution in a subspace $\widetilde{W}\subset W$ which satisfies a certain number of continuity constraints. We have
\begin{equation*}
  \widetilde{W} = W_{\Pi}\oplus \widetilde{W}_{\Delta} = W_{\Pi}\oplus \Pi_{i = 1}^N\widetilde{W}_i,
\end{equation*}
where the function spaces are defined as follows:
\begin{align*}
  W_{\Pi} &= \{\v{u}\in W_0| \lambda_e(\v{u}) = C_{ij},e\subset\Gamma_{ij}\in\mathcal{E}_H, 1\le i\neq j\le N\},   \\
  \widetilde{W}_i &= \{\v{u}_i\in W_i| \sum_{e\in\mathcal{E}_h,e\subset\Gamma_{ij}}\lambda_e(\v{u}_i)\cdot\vert e\vert = 0,\Gamma_{ij}\in\mathcal{E}_H,\Gamma_{ij}\subset\partial\Omega_i \}.
\end{align*}
Here $C_{ij}$ is a constant for every $\Gamma_{ij}\in\mathcal{E}_H$.\\
\indent
Then we have
\begin{equation}\label{matNew}
  \begin{pmatrix}
    A_{II} & A_{I\widetilde{\Delta}} & A_{I\Pi} \\
    A_{\widetilde{\Delta} I} & A_{\widetilde{\Delta}\widetilde{\Delta}}+\gamma M_{\widetilde{\Delta}\widetilde{\Delta}} & A_{\widetilde{\Delta}\Pi} \\
    A_{\Pi I} & A_{\Pi\widetilde{\Delta}} & A_{\Pi\Pi}+\gamma M_{\Pi\Pi} 
  \end{pmatrix}
  \begin{pmatrix}
    u_I \\
    u_{\widetilde{\Delta}} \\
    u_{\Pi} 
  \end{pmatrix} = 
  \begin{pmatrix}
    f_I \\
    f_{\widetilde{\Delta}}+M_{\widetilde{\Delta}\widetilde{\Delta}}g_{\widetilde{\Delta}} \\
    f_{\Pi}+M_{\Pi\Pi}g_{\Pi} 
  \end{pmatrix},
\end{equation}
where $\v{u}_{\Pi}\in W_{\Pi}$, $\v{u}_I,\v{u}_{\widetilde{\Delta}}\in \widetilde{W}_{\Delta}$ and the subscripts $I,\widetilde{\Delta}$ denote the interior degrees of freedom and those on $\Gamma$. Besides, the Robin boundary condition $g$ is also decomposed into two parts, $g_{\widetilde{\Delta}}$ and $g_{\Pi}$. \\
\indent
Now we could describe our new algorithm as follows:
\begin{definition}[Iterative algorithm for H(div)-elliptic problem]\label{defAlgo}
  Given $g^0 = \begin{pmatrix}
                 g_{\widetilde{\Delta}}^0 \\
                 g_{\Pi}^0 
               \end{pmatrix}( = 0)\in V$, we perform the following steps until convergence:\\
  \textbf{Step 1} Solve the Robin boundary problem
  \begin{equation*}
    \begin{pmatrix}
    A_{II} & A_{I\widetilde{\Delta}} & A_{I\Pi} \\
    A_{\widetilde{\Delta} I} & A_{\widetilde{\Delta}\widetilde{\Delta}}+\gamma M_{\widetilde{\Delta}\widetilde{\Delta}} & A_{\widetilde{\Delta}\Pi} \\
    A_{\Pi I} & A_{\Pi\widetilde{\Delta}} & A_{\Pi\Pi}+\gamma M_{\Pi\Pi} 
  \end{pmatrix}
  \begin{pmatrix}
    u_I^n \\
    u_{\widetilde{\Delta}}^n \\
    u_{\Pi}^n 
  \end{pmatrix} = 
  \begin{pmatrix}
    f_I \\
    f_{\widetilde{\Delta}}+M_{\widetilde{\Delta}\widetilde{\Delta}}g_{\widetilde{\Delta}}^n \\
    f_{\Pi}+M_{\Pi\Pi}g_{\Pi}^n 
  \end{pmatrix},
  \end{equation*}
  \noindent
  \textbf{Step 2} Update the Robin boundary condition by equality
  \begin{equation}\label{update}
    \begin{pmatrix}
      \tilde{g}^n_{ij,\widetilde{\Delta}} \\
      \tilde{g}^n_{ji,\widetilde{\Delta}} 
    \end{pmatrix} = 2\gamma
    \begin{pmatrix}
      u_{ji,\widetilde{\Delta}}^n \\
      u_{ij,\widetilde{\Delta}}^n 
    \end{pmatrix}-
    \begin{pmatrix}
      g_{ij,\widetilde{\Delta}}^n \\
      g_{ji,\widetilde{\Delta}}^n 
    \end{pmatrix},
  \end{equation}
  \noindent
  \textbf{Step 3} Get the next iterate by a relaxation
  \begin{equation*}
    g^{n+1} = \theta \tilde{g}^n+(1-\theta)g^n.
  \end{equation*}
\end{definition}

\subsection{Implementation of the algorithm}
There are two ways to implement our new algorithm. The first way is to express functions in $W_{\Pi}$ and $\widetilde{W}_{\Delta}$ explicitly. By the definition of the function spaces, we could regard $u_{\Pi}$ as the edge average degree of freedom and $u_{\Delta}$ as the residual part. Therefore we can change variables on $\Gamma$ to make them explicit. We refer the readers interested in the explicit way to \cite{li2006Cholesky} for further details.\\
\indent
The second way is to introduce a Lagrange multiplier $\mu^n$ to enforce the continuity of edge average part. Then the linear system (\ref{matNew}) could be written in the following format:
\begin{equation}\label{matLag}
  \begin{pmatrix}
    A_{II} & A_{I\Delta} & 0 \\
    A_{\Delta I} & A_{\Delta\Delta}+\gamma M & B^{T} \\
    0 & B & 0 
  \end{pmatrix}
  \begin{pmatrix}
    u_I^n \\
    u_{\Delta}^n \\
    \mu^n 
  \end{pmatrix} = 
  \begin{pmatrix}
    f_I \\
    f_{\Delta}+Mg^n \\
    0 
  \end{pmatrix},
\end{equation}
where
\begin{equation*}
  Bu_{\Delta}^n|_{\Gamma_{ij}} = \int_{\Gamma_{ij}}(\v{u}_{ij,\Delta}^n-\v{u}_{ji,\Delta}^n)\cdot\v{n}_{ij}\text{d}s.
\end{equation*}

To solve the linear system (\ref{matLag}), we first eliminate variables $u_I^n, u_{\Delta}^n$ and get
\begin{equation*}
  S_{\Pi}\mu^n = f_{\Pi}^n
\end{equation*}
where
\begin{gather*}
  S_{\Pi} = 
  \begin{pmatrix}
    0 \\
    B^T 
  \end{pmatrix}^T
  \begin{pmatrix}
    A_{II} & A_{I\Delta} \\
    A_{\Delta I} & A_{\Delta\Delta}+\gamma M 
  \end{pmatrix}^{-1}
  \begin{pmatrix}
    0 \\
    B^T 
  \end{pmatrix},\\
  f_{\Pi}^n = 
  \begin{pmatrix}
    0 \\
    B^T 
  \end{pmatrix}^T
  \begin{pmatrix}
    A_{II} & A_{I\Delta} \\
    A_{\Delta I} & A_{\Delta\Delta}+\gamma M 
  \end{pmatrix}^{-1}
  \begin{pmatrix}
    f_I \\
    f_{\Delta}+Mg^n 
  \end{pmatrix}.
\end{gather*}
We note that the $2\times 2$ block matrix that is inversed can be made block diagonal, with each block associated with a Robin boundary problem in a subdomain. So we can do the elimination at a modest cost. The resulting system of Lagrange multiplier $\mu^n$ is small and it could be solved directly. Then we substitute $\mu$ into (\ref{matLag}) and solve local problems again to obtain the desired vector $u_{\Delta}^n$.
\begin{remark}\label{remark1}
  We could get the operator $Q$ of the iteration as follows: First, eliminate variables $u_I$ and get
  \begin{equation*}
    \begin{pmatrix}
      S+\gamma M & B^T \\
      B & 0 
    \end{pmatrix}
    \begin{pmatrix}
      u_{\Delta}^n \\
      \mu^n 
    \end{pmatrix} = 
    \begin{pmatrix}
      \tilde{f}+Mg^n \\
      0 
    \end{pmatrix}
  \end{equation*}
  where
  \begin{gather*}
    S = A_{\Delta\Delta}-A_{\Delta I}A_{II}^{-1}A_{I\Delta},\\
    \tilde{f} = f_{\Delta}-A_{\Delta I}A_{II}^{-1}f_I.
  \end{gather*}
  We next compute the Lagrange multiplier by the following system
  \begin{equation*}
    B(S+\gamma M)^{-1}B^{T}\mu^n = B(S+\gamma M)^{-1}(\tilde{f}+Mg^n).
  \end{equation*}
  Then the desired vector $u_{\Delta}^n$ could be expressed as
  \begin{equation}\label{udelta}
    u_{\Delta}^n = [(S+\gamma M)^{-1}-K](\tilde{f}+Mg),
  \end{equation}
  where 
  \begin{equation*}
    K = (S+\gamma M)^{-1}B^T[B(S+\gamma M)B^{-1}]^{-1}B(S+\gamma M)^{-1}.
  \end{equation*}
  Applying the equality (\ref{update}) and the relaxation step, omit the same quantity in every iteration, we may get
  \begin{equation*}
    Q = (1-\theta)I+\theta T[2\gamma(S_M^{-1}-K)M-I],
  \end{equation*}
  where the operator $T$ acts on $u_{\Delta}$ as follows:
  \begin{equation*}
    T\begin{pmatrix}
       u_{ij,\Delta} \\
       u_{ji,\Delta} 
     \end{pmatrix} = 
     \begin{pmatrix}
       u_{ji,\Delta} \\
       u_{ij,\Delta} 
     \end{pmatrix}
  \end{equation*}
  and we abbreviate $S+\gamma M$ to $S_M$ for simplicity.\\
  \indent
  In the next section, we will compute the eigenvalues of operator $Q$ to show the optimality of our algorithm.
\end{remark}
\begin{remark}\label{remark2}
  Recall the equality used in step 2, we find that for the numerical solution $u$ and Robin boundary condition $g$ on $\Gamma$ of original problem, it holds that
  \begin{equation}\label{bdeq}
    Tg = 2\gamma u_{\Delta}-g.
  \end{equation}
  If we substitute the result (\ref{udelta}) into (\ref{bdeq}), we may get the following linear system
  \begin{equation*}
    Gg = f_g,
  \end{equation*}
  where
  \begin{gather*}
    G = (TM+M)-2\gamma M(S_M^{-1}-K)M \\
    f_g = 2\gamma M(S_M^{-1}-K)\tilde{f}.
  \end{gather*}
  We will solve the symmetric system of Robin boundary condition $g$ by MINRES in the next section. Besides, we also test the algorithm on Poisson equation and Stokes equations with linear systems derived in the same way.
\end{remark}

\section{Numerical experiments}
In this section, we carry out some numerical experiments to test our algorithm.\\
\indent
We solve the model equation on the unit square $\Omega = [0,1]^2$. The domain $\Omega$ is divided into $N\times N$ rectangular subdomains. The coefficient $\beta$ is set to be $1$. The exact solution is chosen as
\begin{equation*}
\v{u}(x,y) = 
  \begin{pmatrix}
    x(1-x) \\
    y(1-y) 
  \end{pmatrix}.
\end{equation*}

\indent
We first test the iterative algorithm with fixed ratio $H/h = 8$ or fixed number of subdomains. The stopping criterion for the iteration is
\begin{equation*}
  \Vert g^n-g^{n-1}\Vert_{l^{\infty}} < 10^{-6}.
\end{equation*}

\begin{table}[h]
\caption{The iteration numbers and numerical errors with fixed $\frac Hh = 8$}
  \centering\label{tableratio8}
    \begin{tabular}{ccccccc}
    \toprule
    \toprule
$N\times N$& $\gamma = h,\theta = \frac 12$ & $\gamma = h,\theta = \frac 23$ &$\gamma = H,\theta = \frac 12$ & $\gamma = H,\theta = \frac 23$ &$L^2$-err&H(div)-err \cr
    \midrule
$4\times 4$&	41& 30 & 73 & 54 &6.160e-3&1.804e-2 \cr
$8\times 8$&	39&	28 & 68 & 51 &3.081e-3&9.020e-3 \cr
$16\times 16$&	36&	26 & 62 & 46 &1.541e-3&4.510e-3 \cr
$24\times 24$&  35& 25 & 58 & 43 &1.027e-3&3.007e-3 \cr
$32\times 32$&	34& 24 & 55 & 41 &7.706e-4&2.255e-3 \cr
$40\times 40$&  33& 24 & 53 & 40 &6.168e-4&1.804e-3 \cr
$48\times 48$&  33& 23 & 52 & 38 &5.141e-4&1.504e-3 \cr
$64\times 64$&  32& 23 & 49 & 37 &3.865e-4&1.128e-3 \cr
    \bottomrule
    \bottomrule
    \end{tabular}
\end{table}

\begin{table}[h]
\caption{The iteration numbers with fixed $4\times 4$ subdomains}
  \centering\label{tabledomains44}
    \begin{tabular}{ccccc}
    \toprule
    \toprule
$H/h$& $\gamma = h,\theta = \frac 12$ & $\gamma = h,\theta = \frac 23$ &$\gamma = H,\theta = \frac 12$ & $\gamma = H,\theta = \frac 23$  \cr
    \midrule
$4$&	26& 18 & 45 & 33  \cr
$8$&	41&	30 & 73 & 54  \cr
$16$&	72&	53 & 115 & 86  \cr
$32$&	128& 95 & 178 & 133  \cr
    \bottomrule
    \bottomrule
    \end{tabular}
\end{table}

The numerical results presented in Table \ref{tableratio8} and Table \ref{tabledomains44} provide the following key insights. First, from the $L^2$-err and $H(\text{div})$-err data in Table \ref{tableratio8}, we may see that the numerical results converge to true solution as the mesh size $h$ decreases and the reduction rates of errors satisfy the standard a priori estimates of Raviart-Thomas finite element approximation. This validates the correctness of our algorithm. Second, the iteration numbers with different $\gamma$ decrease gradually and converge to a stable number asymptotically as the number of subdomains increases. While, if we fix the number of subdomains and increase the ratio between diameter of subdomains and mesh size, the iteration numbers increase rapidly. As a numerical conclusion, the convergence rate depends solely on the ratio $H/h$. Besides, the choices of relaxation parameter $\theta$ and Robin parameter $\gamma$ significantly influence the convergence rate of the algorithm.

As is demonstrated in Remark \ref{remark1}, the iteration operator $Q$ could be expressed explicitly. To better demonstrate the algorithm's optimality, the following figures present the eigenvalue distribution of the operator, revealing its spectral characteristics which underlie the convergence behavior.

\begin{figure}[H]
    \centering
    \subfigure[$H/h = 4$]{
        \includegraphics[width=0.48\linewidth]{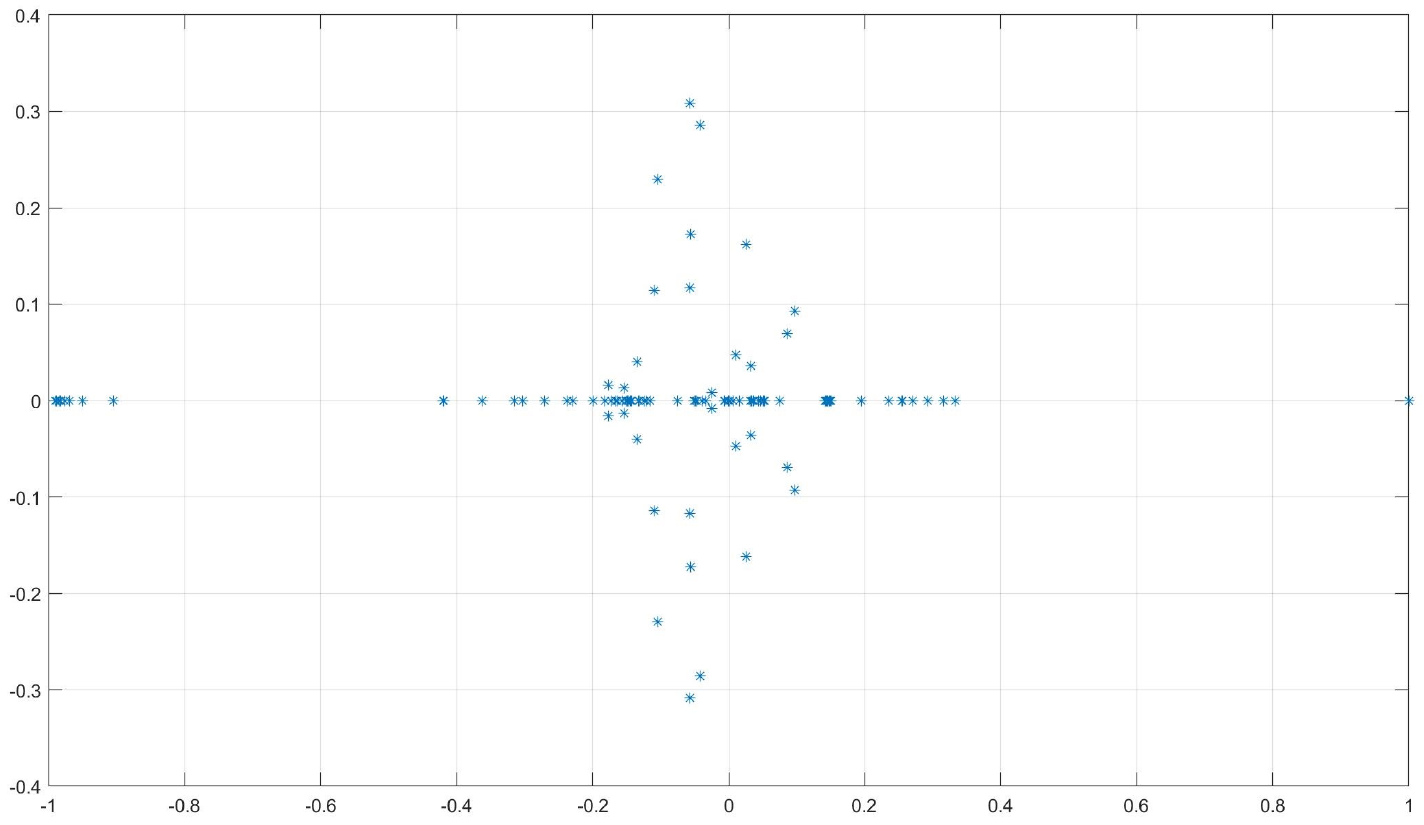}
    }
    \centering
    \subfigure[$H/h = 8$]{
        \includegraphics[width=0.48\linewidth]{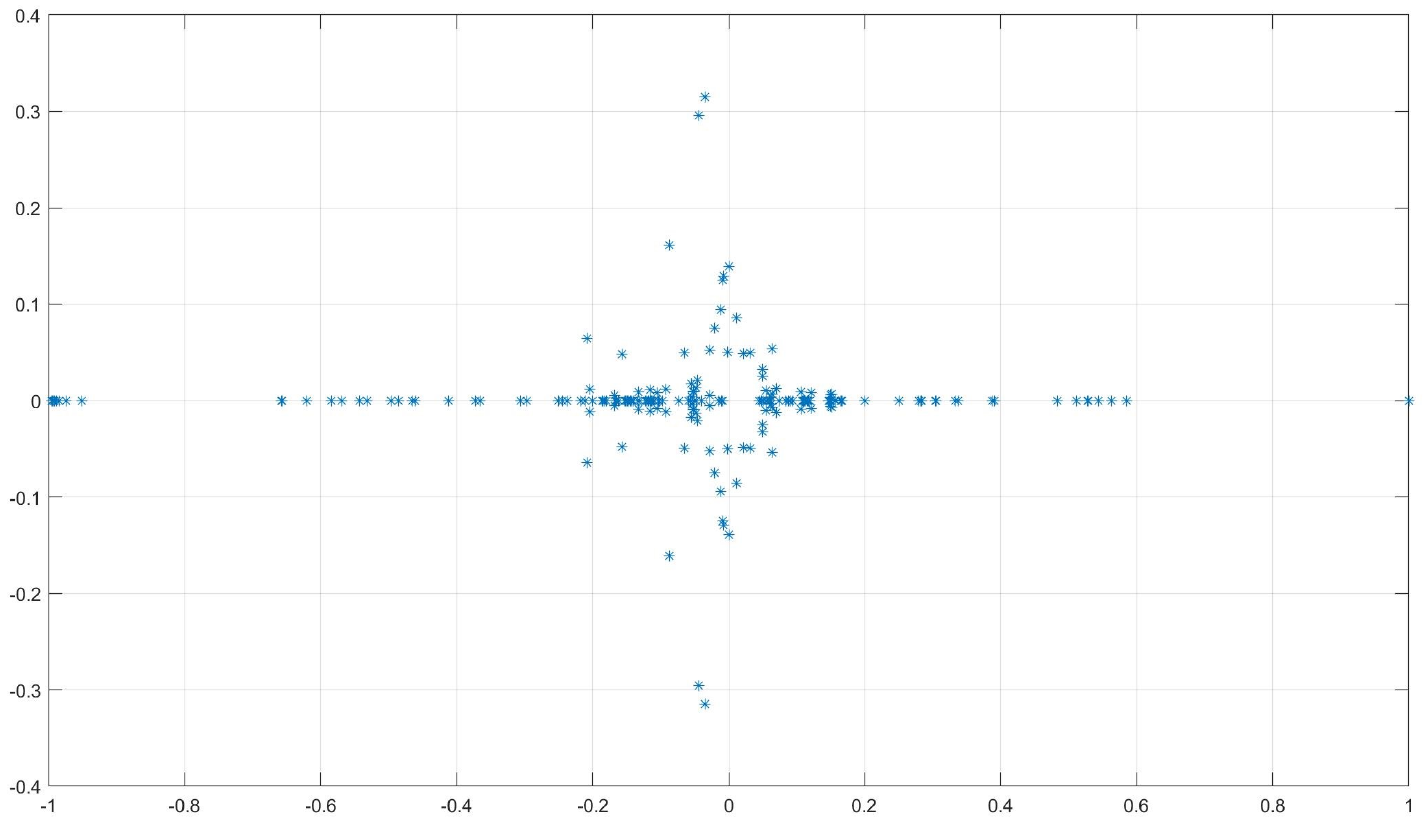}
    }
    \caption{$4\times 4$ subdomains}
    \label{fig1}
\end{figure}

\begin{figure}[H]
    \centering
    \subfigure[$H/h = 4$]{
        \includegraphics[width=0.48\linewidth]{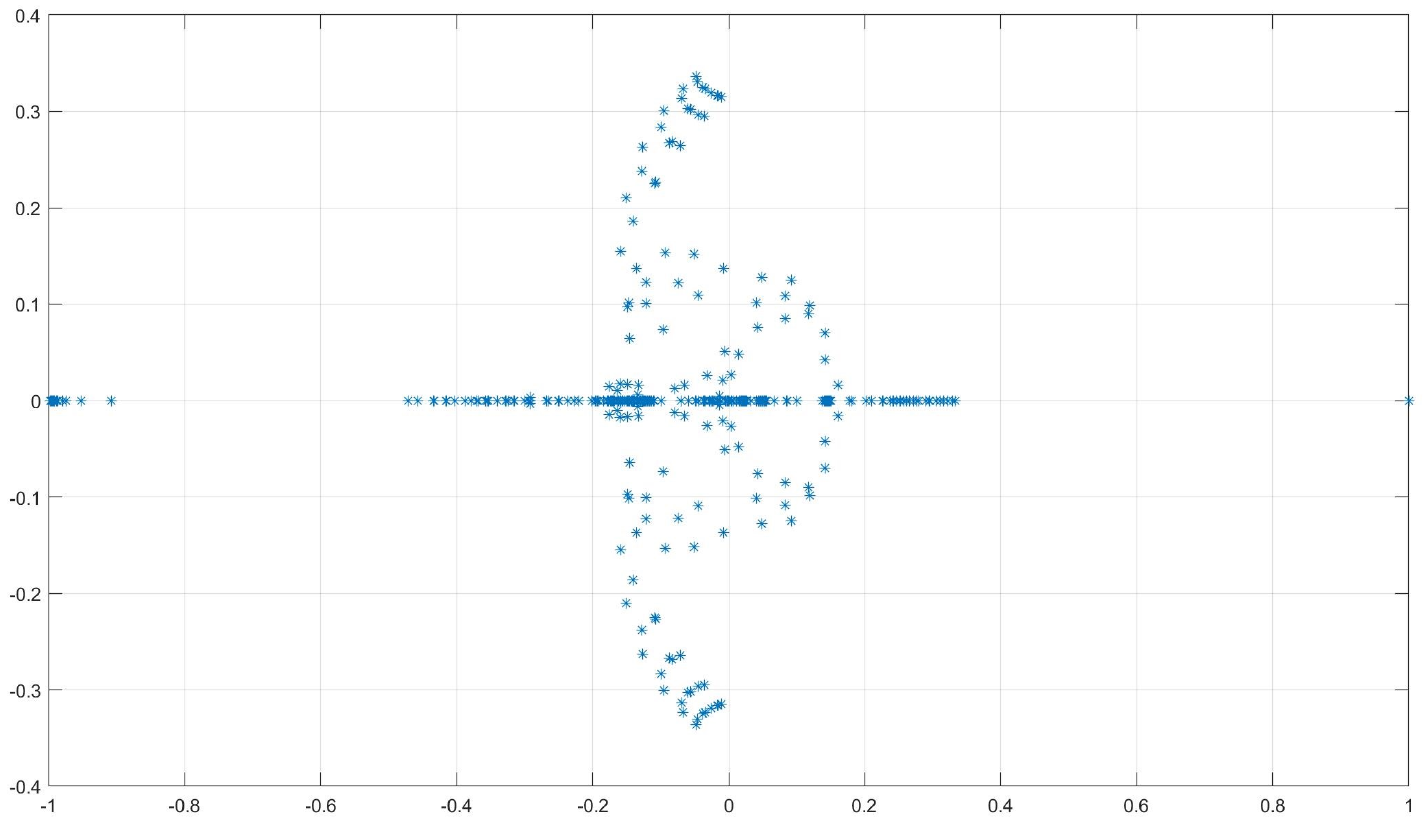}
    }
    \centering
    \subfigure[$H/h = 8$]{
        \includegraphics[width=0.48\linewidth]{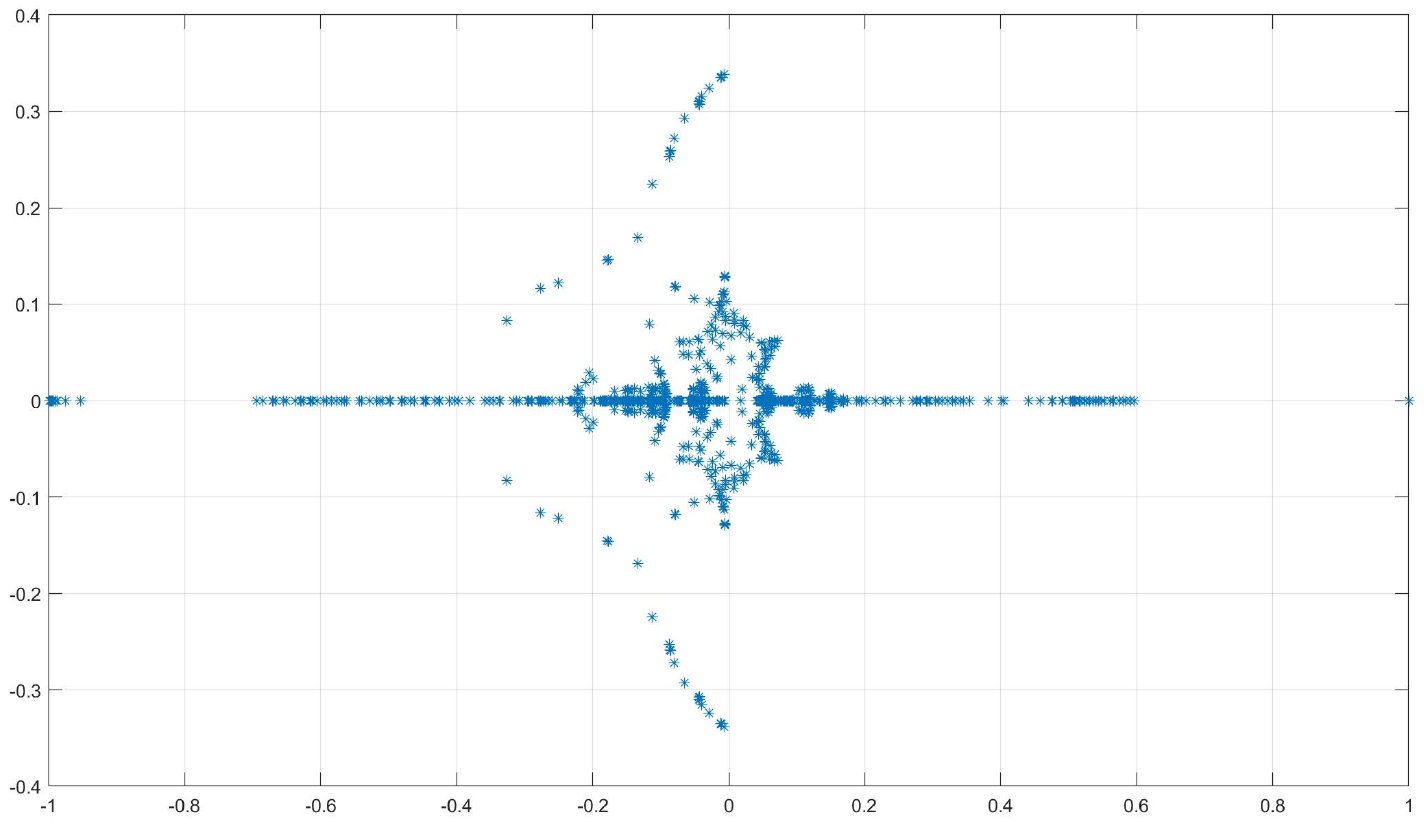}
    }
    \caption{$8\times 8$ subdomains}
    \label{fig2}
\end{figure}

\begin{figure}[H]
    \centering
    \subfigure[$H/h = 4$]{
        \includegraphics[width=0.48\linewidth]{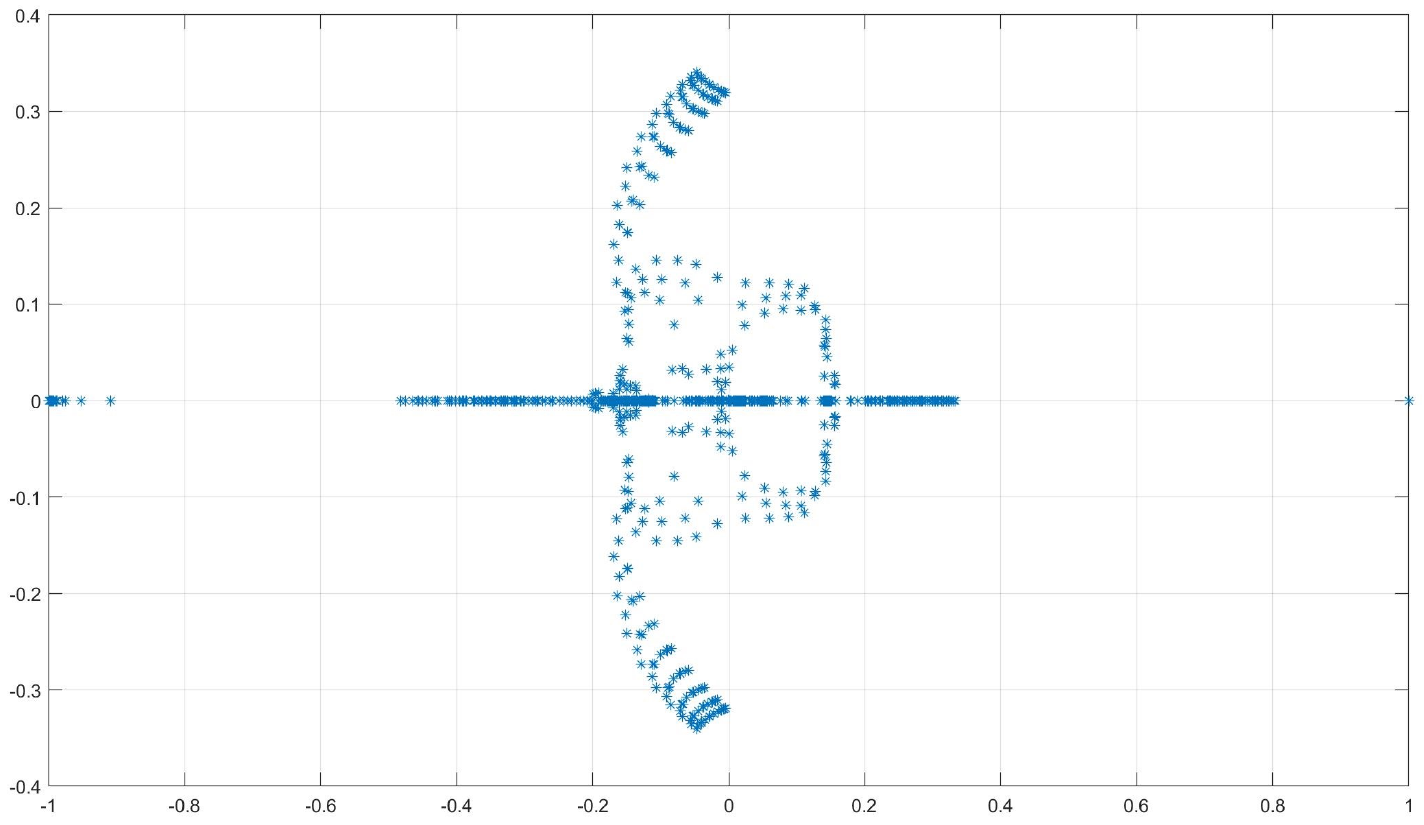}
    }
    \centering
    \subfigure[$H/h = 8$]{
        \includegraphics[width=0.48\linewidth]{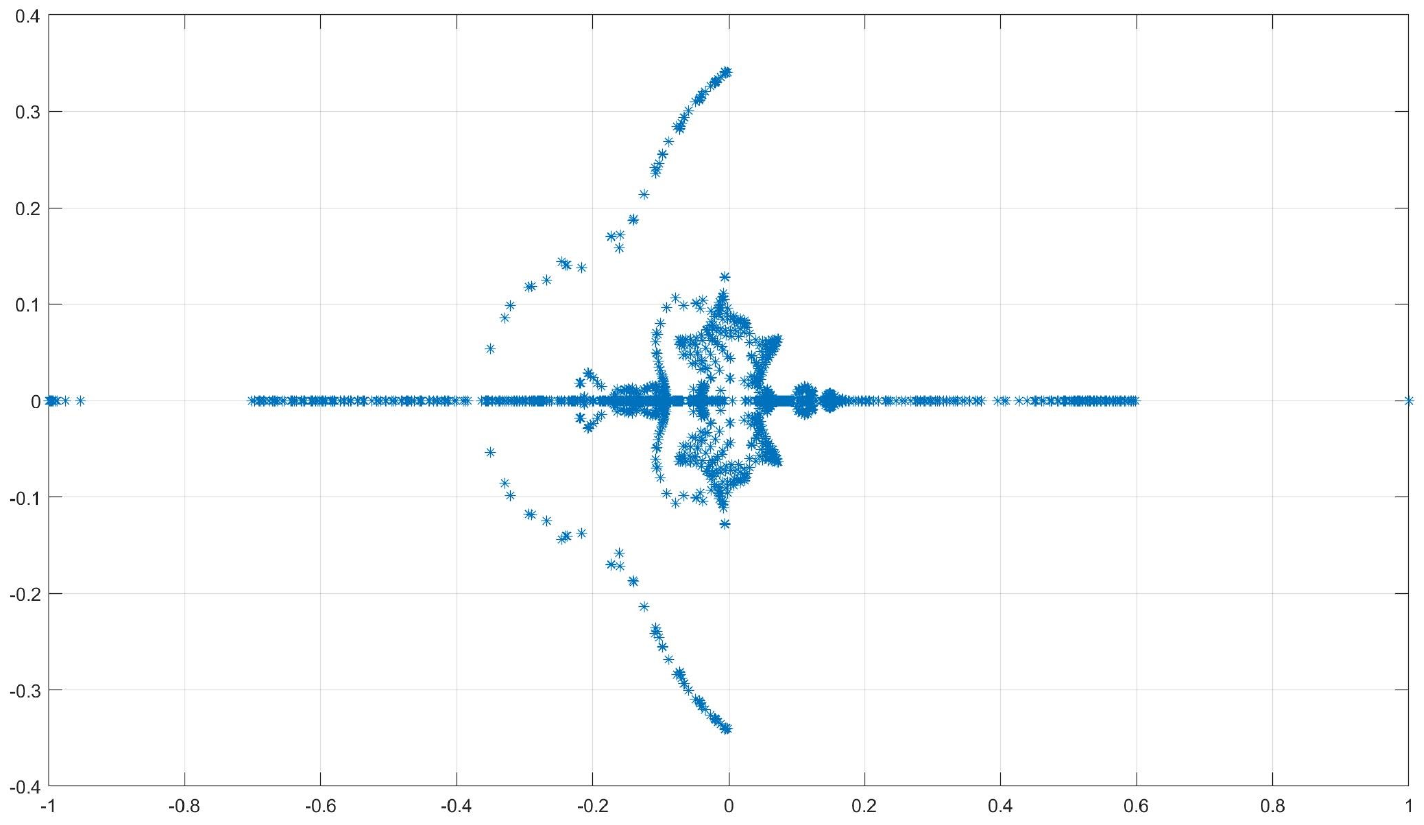}
    }
    \caption{$12\times 12$ subdomains}
    \label{fig3}
\end{figure}

Here we compute the eigenvalues of iteration operator without relaxation for simplicity. From the figures above, we may get some conclusions:
\begin{enumerate}
  \item The iteration operator has an eigenvalue equal to 1 with arbitrary number of subdomains and ratio between diameter of subdomains $H$ and mesh size $h$. The eigenvalue corresponds to the true solution.
  \item The complex eigenvalues exhibit bounded moduli.
  \item The upper bound of real eigenvalues except 1 depends on the ratio $H/h$ while it is independent of the number of subdomains.
\end{enumerate}
The numerical evidence additionally confirms that if we fix the ratio $H/h$, the moduli of the eigenvalues with different number of subdomains have a common bound by applying an appropriate relaxation.

\begin{table}[H]
\caption{The iteration number of MINRES with different $\gamma$ and $H/h$}
  \centering\label{table2}
    \begin{tabular}{ccccc}
    \toprule
    \toprule
$N\times N$& $\gamma = h, H/h = 8$ &$\gamma = h, H/h = 16$&$\gamma = H, H/h = 8$&$\gamma = H, H/h = 16$\cr
    \midrule
$4\times 4$&	31 &45&30&39\cr
$8\times 8$&	32 &45&34&43\cr
$16\times 16$&	31 &45&32&44\cr
$24\times 24$&  30 &43&32&44\cr
$32\times 32$&	30 &43&32&41\cr
$40\times 40$&  30 &41&32&41\cr
$48\times 48$&  30 &41&32&41\cr
    \bottomrule
    \bottomrule
    \end{tabular}
\end{table}

\begin{table}[H]
\caption{The iteration number of MINRES (Poisson and Stokes equations)}
  \centering\label{table3}
    \begin{tabular}{ccc}
    \toprule
    \toprule
$N\times N$& Poisson & Stokes \cr
    \midrule
$4\times 4$&	16 &29\cr
$8\times 8$&	23 &31\cr
$16\times 16$&	23 &28\cr
$24\times 24$&  21 &24\cr
$32\times 32$&	19 &22\cr
$40\times 40$&  17 &20\cr
$48\times 48$&  18 &21\cr
    \bottomrule
    \bottomrule
    \end{tabular}
\end{table}

Finally, we solve the linear system mentioned in Remark \ref{remark2} by MINRES. The MINRES iteration stops when $\Vert Gg^n-f_g\Vert_{l^2}/\Vert f_g\Vert_{l^2}$ is less than $10^{-6}$. The data presented in Table \ref{table2} show that for fixed ratio $H/h$, the iteration counts of MINRES method approach a stable number, which is similar with the behaviour of previous iterative method. Besides, the iteration numbers in Table \ref{table3} confirm that the algorithm solved by MINRES remains effective when applied to Poisson equation and Stokes equations.

\section{Conclusion}
In this paper, we develop a new two-side Robin-Robin domain decomposition method for H(div)-elliptic method and derive an algebraic system for Robin boundary conditions from the iterative algorithm. Numerical results show that the convergence rate of iterative algorithm and the condition number of the algebraic system depend solely on $H/h$. Besides, we could get similar algebraic systems of Poisson equation and Stokes equations and both of them have condition numbers only depend on $H/h$. There still remains some problems, such as the relationship between convergence rate and the ratio $H/h$, how the Robin parameters and relaxation parameter influence the convergence rate and why the algebraic system remains optimality when applied on Poisson equation and Stokes equations while the iterative method fails. These problems will be our future work.

\newpage

\bibliographystyle{plain}
\bibliography{reference}
\end{document}